\documentclass[11pt]{article}

\usepackage{latexsym}
\usepackage{amssymb}

\newtheorem{definition}{Definition}

\newtheorem{remark}{Remark}
\newtheorem{example}{Example}

\begin{document}
{
\begin{center}
{\Large\bf
The operator approach to the truncated multidimensional moment problem.}
\end{center}
%
\section{Introduction.}

Let us introduce some notations. As usual, we denote by $\mathbb{R}, \mathbb{C}, \mathbb{N}, \mathbb{Z}, \mathbb{Z}_+$
the sets of real numbers, complex numbers, positive integers, integers and non-negative integers,
respectively.
By $\mathbb{Z}^n_+$ we mean $\mathbb{Z}_+\times \ldots \times\mathbb{Z}_+$, and $\mathbb{R}^n = \mathbb{R}\times \ldots \times\mathbb{R}$,
where the decart products are taken with $n$ copies.
Let $\mathbf{k} = (k_1,\ldots,k_n)\in\mathbb{Z}^n_+$, $\mathbf{t} = (t_1,\ldots,t_n)\in\mathbb{R}^n$. Then
$\mathbf{t}^{\mathbf{k}}$ means the monomial $t_1^{k_1}\ldots t_n^{k_n}$, and $|\mathbf{k}| = k_1 + \ldots + k_n$.
By $\mathfrak{B}(\mathbb{R}^n)$ we denote the set of all Borel subsets of $\mathbb{R}^n$.

Let $\mathcal{K}$ be an arbitrary finite subset of $\mathbb{Z}^n_+$. Let $\mathcal{S} = (s_{\mathbf{k}})_{\mathbf{k}\in\mathcal{K}}$
be an arbitrary set of real numbers.
\textit{The truncated multidimensional moment problem} consists of finding a (non-negative) measure $\mu$ on $\mathfrak{B}(\mathbb{R}^n)$
such that
\begin{equation}
\label{f1_1}
\int \mathbf{t}^{\mathbf{k}} d\mu(\mathbf{t}) = s_{\mathbf{k}},\qquad  \forall \mathbf{k}\in\mathcal{K}.
\end{equation}

The multidimensional moment problem (both the full and the truncated versions) turned out to be much more complicated than its
one-dimensional prototype~\cite{cit_1000_Berezansky_1965__Book}, \cite{cit_980_Berg_Christiansen_Ressel__Book}, \cite{cit_990_Marshall_Book}. 
An operator-theoretical interpretation of the (full) 
multidimensional moment problem was given
by Fuglede in~\cite{cit_980_Fuglede}. It should be noticed that the operator approach to moment problems was introduced
by Naimark in 1940--1943 and then developed by many authors, see historical notes in~\cite{cit_15000_Zagorodnyuk_2017_J_Adv_Math_Stud}.
Elegant conditions for the solvability of the multidimensional moment problem in the case of the support on semi-algebraic sets were given by
Schm\"udgen in~\cite{cit_995_Schmudgen__1991}, \cite{cit_995_Schmudgen__2003}.
Another conditions for the solvability of the multidimensional moment problem, using an extension of the moment sequence, 
were given by Putinar and Vasilescu, see~\cite{cit_993_P_V__1999}, \cite{cit_14100_Vasilescu}.
Developing the idea of Putinar and Vasilescu, we presented different conditions for the solvability of
the two-dimensional moment problem and proposed an algorithm (which essentially consists of solving of linear equations)
for a construction of the solutions set~\cite{cit_14500_Zagorodnyuk_2010_AFA}.
An analytic parametrization for all solutions of the two-dimensional moment problem in a strip was given
in~\cite{cit_14700_Zagorodnyuk_2013_MFAT}. 
Another approach to multidimensional and complex moment problems (including truncated problems), using extension arguments for $*$-semigroups, 
has been developed by Cicho\'n, Stochel and Szafraniec, see~\cite{cit_982_C_St_Sz_2011} and references therein.
Still another approach for the two-dimensional moment problem was proposed by Ovcharenko
in~\cite{cit_991_Ovcharenko_1}, \cite{cit_991_Ovcharenko_2}.

In this paper we shall be focused on the truncated multidimensional moment problem.
A general approach for this moment problem was given by Curto and Fialkow in their books~\cite{cit_985_Curto_Fialkow__Book1} 
and~\cite{cit_985_Curto_Fialkow__Book2}.
These books entailed a series of papers by a group of mathematicians, see recent papers~\cite{cit_987_Fialkow_2011}, 
\cite{cit_14200_Vasilescu}, \cite{cit_14000_Yoo}
and references therein.
This approach includes an extension of the matrix of prescribed moments with the same rank. 
Effective optimization algorithms for the multidimensional moment problems were given in the book of Lasserre~\cite{cit_985_Lasserre_Book}.
Another approach for truncated moment problems, using a notion of an idempotent, was presented by Vasilescu in~\cite{cit_14220_Vasilescu}.
Atomic solutions to various matrix truncated $K$-moment problems were studied by Kimsey and Woerdeman in~\cite{cit_983_Kimsey_Woerdeman}.
There exists a connection of the truncated multidimensional moment problems with the completion problems for subnormal
operators, see, e.g.,~\cite{cit_984_Kimsey}. Observe that the complexification of the real truncated moment problem needs the even dimension $d$.
The complexification of the truncated multidimensional moment problem and the use of hyponormal operators was investigated
by Kimsey and Putinar
in~\cite{cit_982_Kimsey_Putinar}.
We should also mention recent papers~\cite{cit_15000_Zagorodnyuk_2018_AOT}, \cite{cit_998_Dio_Schmudgen_ArXiv} 
on the subject.

Even in the one-dimensional case ($n=1$) the operator approach is effective not for all types of truncations $\mathcal{K}$.
Thus, we need to define some admissible types of truncations, where we can get solutions. The second feature of the truncated case is that
we need to take care that the corresponding multiplication operators in the associated Hilbert space are well-defined.
If all the above is done, we come to a problem of an extension of commuting symmetric operators.
In the case of the dimensional stability (see a precise definition below), the corresponding operators are self-adjoint and we
get an atomic solution.
More weak conditions which give a way for an explicit check for solutions are given, as well.

\noindent
{\bf Notations. }
Besides the given above notations we shall use the following conventions.
By $\mathbb{Z}_{k,l}$ we mean all integers $j$ satisfying the following inequality:
$k\leq j\leq l$; ($k,l\in\mathbb{Z}$).
If H is a Hilbert space then $(\cdot,\cdot)_H$ and $\| \cdot \|_H$ mean
the scalar product and the norm in $H$, respectively.
Indices may be omitted in obvious cases.
For a linear operator $A$ in $H$, we denote by $D(A)$
its  domain, by $R(A)$ its range, and $A^*$ means the adjoint operator
if it exists. If $A$ is invertible then $A^{-1}$ means its
inverse. $\overline{A}$ means the closure of the operator, if the
operator is closable. If $A$ is bounded then $\| A \|$ denotes its
norm.
For a set $M\subseteq H$
we denote by $\overline{M}$ the closure of $M$ in the norm of $H$.
By $\mathop{\rm Lin}\nolimits M$ we mean
the set of all linear combinations of elements from $M$,
and $\mathop{\rm span}\nolimits M:= \overline{ \mathop{\rm Lin}\nolimits M }$.
By $E_H$ we denote the identity operator in $H$, i.e. $E_H x = x$,
$x\in H$. In obvious cases we may omit the index $H$. If $H_1$ is a subspace of $H$, then $P_{H_1} =
P_{H_1}^{H}$ is an operator of the orthogonal projection on $H_1$
in $H$.

\section{Necessary conditions for the solvability of the moment problem.}

Consider the following operator $W_j$ on $\mathbf{Z}^n_+$:
\begin{equation}
\label{f2_1}
W_j (k_1, \ldots, k_{j-1}, k_j, k_{j+1},\ldots, k_n) = (k_1, \ldots, k_{j-1}, k_j + 1, k_{j+1},\ldots, k_n),
\end{equation}
for $j=1,\ldots,n$. Thus, the operator $W_j$ increases the $j$-th coordinate.

Probably, the following kind of subsets appeared for the first time in the work of Kimsey and Woerdeman~\cite{cit_983_Kimsey_Woerdeman}.
\begin{definition}
\label{d2_1}
A finite subset $K\subset \mathbb{Z}^n_+$ is said to be \textbf{admissible}, if the following conditions hold:

\begin{itemize}
\item[1)] $\mathbf{0} = (0,\ldots,0)\in K$;

\item[2)] $\forall \mathbf{k}\in K\backslash\{ \mathbf{0} \}$, 
\begin{equation}
\label{f2_5}
\mathbf{k} = W_{a_{|\mathbf{k}|}} W_{a_{|\mathbf{k}| - 1}} \ldots W_{a_1} \mathbf{0},
\end{equation}
for some $a_j\in\{ 1,\ldots,n \}$, and
\begin{equation}
\label{f2_7}
\widetilde{\mathbf{k}}_r := W_{a_r} \ldots W_{a_1} \mathbf{0} \in K,\qquad \forall r=1,2,\ldots, |\mathbf{k}|.
\end{equation}
\end{itemize}

\end{definition}

\begin{example}
\label{e1_1}
1) Let $K = K_r = \{ \mathbf{k}\in\mathbb{Z}^n_+: |\mathbf{k}|\leq r \}$, $r\in\mathbb{Z}_+$. Then $K$ is admissible, since
$\forall \mathbf{k} = (k_1,\ldots,k_n) \in K\backslash\{ \mathbf{0} \}$,
$$ \mathbf{k} = W_n^{k_n} W_{n-1}^{k_{n-1}} \ldots W_1^{k_1} \mathbf{0}. $$

\noindent
2) Let $K = K_{d_1,d_2,\ldots,d_n} = \{ \mathbf{k}=(k_1,\ldots,k_n)\in\mathbb{Z}^n_+: k_1\leq d_1,\ldots,k_n\leq d_n \}$, 
$d_1,\ldots, d_n\in\mathbb{Z}_+$. 
Notice that the truncated two-dimensional moment problem with rectangular data appeared in~\cite{cit_982_Kimsey_Putinar},
\cite{cit_15000_Zagorodnyuk_2018_AOT}.
The general case of the set $K_{d_1,d_2,\ldots,d_n}$ was proposed to the author by Prof. Vasilescu.

The set $K_{d_1,d_2,\ldots,d_n}$ is admissible, since
$\forall \mathbf{k} = (k_1,\ldots,k_n) \in K\backslash\{ \mathbf{0} \}$,
$$ \mathbf{k} = W_n^{k_n} W_{n-1}^{k_{n-1}} \ldots W_1^{k_1} \mathbf{0}. $$
\end{example}

Suppose that for an admissible finite set $K\subset\mathbb{Z}^n_+$ the moment problem~(\ref{f1_1}), with
$\mathcal{K} = K+K$ and some $\mathcal{S} = (s_{\mathbf{k}})_{\mathbf{k}\in\mathcal{K}}$,
has a solution $\mu$. Let us investigate which properties of the data $\mathcal{S}$ this fact yields.

The first property is the usual positivity condition. For practical purposes, it is useful to introduce some
indexation in the set $K$ by a unique index $j$.
\begin{equation}
\label{f2_9}
K = \left\{
\mathbf{k}_j(\in\mathbb{Z}^n_+),\quad j=0,1,\ldots,\rho 
\right\}.
\end{equation}
Of course, $\rho + 1$ is the number of elements in $K$.
Consider an arbitrary polynomial of the following form:
\begin{equation}
\label{f2_20}
p(\mathbf{t}) = \sum_{j=0}^\rho \alpha_j \mathbf{t}^{\mathbf{k}_j},\qquad \alpha_j\in\mathbb{C}.
\end{equation}
Then
$$ 0 \leq \int |p|^2 d\mu =
\sum_{j,m=0}^\rho \alpha_j \overline{\alpha_m} s_{\mathbf{k}_j + \mathbf{k}_m}. $$
Denote
\begin{equation}
\label{f2_23}
\Gamma = \left( s_{\mathbf{k}_j + \mathbf{k}_m} \right)_{m,j=0}^\rho.
\end{equation}
We obtain the first necessary condition of the solvability:
\begin{equation}
\label{f2_25}
\Gamma \geq 0.
\end{equation}

We now suppose that for an admissible finite set $K\subset\mathbb{Z}^n_+$ the moment problem~(\ref{f1_1}), with
$\mathcal{K} = K+K$ and some $\mathcal{S} = (s_{\mathbf{k}})_{\mathbf{k}\in\mathcal{K}}$, is given and condition~(\ref{f2_25})
holds (we do not require that the moment problem is solvable).
A set $\mathfrak{L}$ of all polynomials of the form~(\ref{f2_20}) is a linear vector space. Let us consider the following functional:
$$ <p, q> = \sum_{j,m=0}^\rho \alpha_j \overline{\beta_m} s_{\mathbf{k}_j + \mathbf{k}_m}, $$
where $p$ is from~(\ref{f2_20}), and $q$ has the same form as $p$, but with $\beta_j (\in\mathbb{C})$ instead of $\alpha_j$.
The functional $<\cdot,\cdot>$ is sesquilinear, $<p,p>\geq 0$, and $\overline{<p,q>} = <q,p>$.
Introducing the classes of the equivalence $[p]_{\mathfrak{L}}$ in $\mathfrak{L}$ we obtain a finite-dimensional Hilbert space $H$.

\textit{We now return to the case of the solvable moment problem}. Consider the space $L^2_\mu$ which consists of (the classes of the equivalence
of) complex-valued measurable functions $f$ such that $\int |f(\mathbf{t})|^2 d\mu < \infty$.
The class of the equivalence in $L^2_\mu$ will be denoted by $[\cdot]_{L^2_\mu}$.
Denote by $T_l$ the multiplication operator in $L^2_\mu$:
\begin{equation}
\label{f2_27}
T_l f(\mathbf{t}) = t_l f(\mathbf{t}),\qquad f\in D_l, 
\end{equation}
where $D_l = \{ f(\mathbf{t})\in L^2_\mu:\ t_l f(\mathbf{t})\in L^2_\mu \}$.

\noindent
Consider the associated Hilbert space $H$, defined as above. The following transformation is useful:
\begin{equation}
\label{f2_29}
W \sum_{j=0}^\rho \alpha_j [ \mathbf{t}^{\mathbf{k}_j} ]_{L^2_\mu} = \sum_{j=0}^\rho \alpha_j [ \mathbf{t}^{\mathbf{k}_j} ]_{\mathfrak{L}},\qquad
\alpha_j\in\mathbb{C}. 
\end{equation}
The transformation $W$ is well-defined, linear and isometric. It maps
$L^2_{\mu;K} := \mathop{\rm Lin}\nolimits \{ [ \mathbf{t}^{\mathbf{k}_j} ]_{L^2_\mu} \}_{j=0}^\rho$
on the whole space $H$. Denote
$\vec e_r := (\delta_{r,m})_{m=1}^n \in \mathbb{Z}^n_+$, $r=1,\ldots,n$, and
\begin{equation}
\label{f2_31}
\Omega_l = \{ j\in \{ 0,\ldots,\rho \}:\ \mathbf{k}_j + \vec e_l \in K \},\qquad l=1,\ldots,n.
\end{equation}
Observe that
$$ W T_l W^{-1} \sum_{j\in\Omega_l} \alpha_j [ \mathbf{t}^{\mathbf{k}_j} ]_{\mathfrak{L}} =
\sum_{j\in\Omega_l} \alpha_j [ \mathbf{t}^{\mathbf{k}_j + \vec e_l} ]_{\mathfrak{L}},\quad \alpha_j\in\mathbb{C}. $$
Since the operator $W T_l W^{-1}$ is well defined, the following implication holds:
\begin{equation}
\label{f2_33}
\left( \sum_{j\in\Omega_l} \alpha_j [ \mathbf{t}^{\mathbf{k}_j} ]_{\mathfrak{L}} = 0,\quad \mbox{for some $\alpha_j\in\mathbb{C}$} \right)
\Rightarrow
\left( \sum_{j\in\Omega_l} \alpha_j [ \mathbf{t}^{\mathbf{k}_j + \vec e_l} ]_{\mathfrak{L}} = 0 \right).
\end{equation}
The latter implication is equivalent to the following one:
$$ \left( \left( \sum_{j\in\Omega_l} \alpha_j [ \mathbf{t}^{\mathbf{k}_j} ]_{\mathfrak{L}}, [ \mathbf{t}^{\mathbf{k}_m} ]_{\mathfrak{L}} \right)_H 
= 0,\ \forall m\in\Omega_l,\quad \mbox{for some $\alpha_j\in\mathbb{C}$} \right)
\Rightarrow $$
\begin{equation}
\label{f2_35}
\left( \left( \sum_{j\in\Omega_l} \alpha_j [ \mathbf{t}^{\mathbf{k}_j + \vec e_l} ]_{\mathfrak{L}}, [ \mathbf{t}^{\mathbf{k}_m + \vec e_l} ]_{\mathfrak{L}} 
\right)_H = 0,\ \forall m\in\Omega_l \right)
\end{equation}
or, equivalently,
$$ \left( \sum_{j\in\Omega_l} \alpha_j s_{\mathbf{k}_j + \mathbf{k}_m} = 0,\ \forall m\in\Omega_l,\quad 
\mbox{for some $\alpha_j\in\mathbb{C}$} \right)
\Rightarrow $$
\begin{equation}
\label{f2_37}
\left( \sum_{j\in\Omega_l} \alpha_j s_{\mathbf{k}_j + \vec e_l + \mathbf{k}_m + \vec e_l}
= 0,\ \forall m\in\Omega_l \right).
\end{equation}
Denote
\begin{equation}
\label{f2_39}
\Gamma_l = \left( s_{\mathbf{k}_j + \mathbf{k}_m} \right)_{m,j\in\Omega_l},\quad
\widehat \Gamma_l = \left( s_{\mathbf{k}_j + \vec e_l + \mathbf{k}_m + \vec e_l} \right)_{m,j\in\Omega_l},\qquad l=1,2,\ldots,n,
\end{equation}
where the indices from $\Omega_l$ are taken in the increasing order.
We obtain the second necessary condition of the solvability:
\begin{equation}
\label{f2_45}
\mathop{\rm Ker}\nolimits \Gamma_l \subseteq \mathop{\rm Ker}\nolimits \widehat \Gamma_l,\qquad l=1,2,\ldots,n.
\end{equation}

\section{The operator approach to the moment problem. The dimensional stability.}

Suppose that for an admissible finite set $K\subset\mathbb{Z}^n_+$ the moment problem~(\ref{f1_1}), with
$\mathcal{K} = K+K$ and some $\mathcal{S} = (s_{\mathbf{k}})_{\mathbf{k}\in\mathcal{K}}$, is given.
Choose and fix some indexation~(\ref{f2_9}). Assume that conditions~(\ref{f2_25}),(\ref{f2_45})
hold.

\noindent
We may construct the associated Hilbert space $H$, as in the previous section. For $l=1,\ldots,n$ we consider the following operators:
\begin{equation}
\label{f3_5}
M_l \sum_{j\in\Omega_l} \alpha_j [ \mathbf{t}^{\mathbf{k}_j} ]_{\mathfrak{L}} =
\sum_{j\in\Omega_l} \alpha_j [ \mathbf{t}^{\mathbf{k}_j + \vec e_l} ]_{\mathfrak{L}},\quad \alpha_j\in\mathbb{C}, 
\end{equation}
with $D(M_l) = \mathop{\rm Lin}\nolimits \{ [ \mathbf{t}^{\mathbf{k}_j} ]_{\mathfrak{L}} \}_{j\in\Omega_l}$.
By condition~(\ref{f2_45}) the operator $M_l$ is well-defined. Moreover, it is linear and symmetric.
In particular, we have
\begin{equation}
\label{f3_7}
M_l [ \mathbf{t}^{\mathbf{k}_j} ]_{\mathfrak{L}} =
[ \mathbf{t}^{\mathbf{k}_j + \vec e_l} ]_{\mathfrak{L}},\qquad j\in\Omega_l;\quad l=1,\ldots,n. 
\end{equation}

Suppose that there exist commuting self-adjoint operators $\widetilde M_j\supseteq M_j$ ($j=1,\ldots,n$) in
a \textit{finite-dimensional} Hilbert space $\widetilde H\supseteq H$. 
Observe that in this case operators $\widetilde M_j$ are bounded and defined on the whole space $\widetilde H$.

Choose an arbitrary $\mathbf{k}=(k_1,\ldots,k_n)\in K\backslash\{ \mathbf{0} \}$. We shall use the notations from Definition~\ref{d2_1}.
By the induction argument one can verify that
\begin{equation}
\label{f3_11}
\left[ \mathbf{t}^{\widetilde{\mathbf{k}}_r} \right]_{\mathfrak{L}} = 
\widetilde M_{a_r} \ldots \widetilde M_{a_1} [1]_{\mathfrak{L}},\qquad r=1,2,\ldots, |\mathbf{k}|.
\end{equation}
In particular, we obtain that
\begin{equation}
\label{f3_15}
\left[ \mathbf{t}^{\mathbf{k}} \right]_{\mathfrak{L}} = 
\widetilde M_{a_{|\mathbf{k}|}} \ldots \widetilde M_{a_1} [1]_{\mathfrak{L}}.
\end{equation}
Since operators $\widetilde M_j$ commute we may rearrange the product in~(\ref{f3_15}). Moreover, it is clear that
$W_1$ appears in the product in~(\ref{f2_5}) $k_1$ times, $W_2$ appears $k_2$ times, ..., $W_n$ appears $k_n$ times.
Then we get
\begin{equation}
\label{f3_17}
\left[ \mathbf{t}^{\mathbf{k}} \right]_{\mathfrak{L}} = 
\widetilde M_{1}^{k_1} \widetilde M_{2}^{k_2} \ldots \widetilde M_{n}^{k_n} [1]_{\mathfrak{L}},\qquad 
\forall \mathbf{k}=(k_1,\ldots,k_n)\in K.
\end{equation}

We can now construct a solution of the moment problem. For an arbitrary $\mathbf{k}=(k_1,\ldots,k_n)\in (K+K)$, 
$\mathbf{k} = \mathbf{k}'+ \mathbf{k}''$, 
$\mathbf{k}'=(k_1',\ldots,k_n'), \mathbf{k}''=(k_1'',\ldots,k_n'')\in K$, we may write
$$ s_{\mathbf{k}} = \left( [\mathbf{t}^{\mathbf{k}'}], [\mathbf{t}^{\mathbf{k}''}] \right)_H = 
\left( \widetilde M_{1}^{k_1'} \ldots \widetilde M_{n}^{k_n'} [1]_{\mathfrak{L}}, 
\widetilde M_{1}^{k_1''} \ldots \widetilde M_{n}^{k_n''} [1]_{\mathfrak{L}} \right)_H = $$
\begin{equation}
\label{f3_18}
= \left( \widetilde M_{1}^{k_1} \ldots \widetilde M_{n}^{k_n} [1]_{\mathfrak{L}}, [1]_{\mathfrak{L}} \right)_H =
\int \mathbf{t}^{\mathbf{k}} d\mu(\mathbf{t}), 
\end{equation}
where
\begin{equation}
\label{f3_19}
\mu(\delta) = \left(
E(\delta) [1]_{\mathfrak{L}}, [1]_{\mathfrak{L}}
\right)_H,\qquad \delta\in\mathfrak{B}(\mathbb{R}^n),
\end{equation}
where $E(\delta)$ is the spectral measure of a commuting tuple $\widetilde M_1,\ldots,\widetilde M_n$.
Consequently, we get a solution $\mu$ of the moment problem.

Denote
\begin{equation}
\label{f3_30}
\Omega_0 = \{ j\in \{0,\ldots,\rho \}:\ \mathbf{k}_j + \vec e_1, \mathbf{k}_j + \vec e_2, \ldots, \mathbf{k}_j + \vec e_n \in K \},
\end{equation}
and
\begin{equation}
\label{f3_35}
H_0 = \mathop{\rm Lin}\nolimits \{ [ \mathbf{t}^{\mathbf{k}_j} ]_{\mathfrak{L}} \}_{j\in\Omega_0}.
\end{equation}
Observe that
$$ \Omega_0\subseteq \Omega_j,\qquad j=1,\ldots,n, $$
and therefore
$$ H_0\subseteq D(M_j),\qquad j=1,\ldots,n. $$

\begin{definition}
\label{d3_1}
Suppose that for an admissible finite set $K\subset\mathbb{Z}^n_+$ the moment problem~(\ref{f1_1}), with
$\mathcal{K} = K+K$ and some $\mathcal{S} = (s_{\mathbf{k}})_{\mathbf{k}\in\mathcal{K}}$, is given and 
conditions~(\ref{f2_25}),(\ref{f2_45})
hold (for some indexation). Define the associated Hilbert space $H$ and its subspace $H_0$.
The set of moments $\mathcal{S}$ is said to be \textbf{dimensionally stable}, if
$\dim H = \dim H_0$.
\end{definition}

Suppose that for the moment problem, as in Definition~\ref{d3_1}, the set $\mathcal{S}$ is dimensionally stable.
Then operators $M_l$ are self-adjoint and defined on the whole $H$. Observe that for $l,r\in\{1,\ldots, n\}:\ l\not= r$, we have 
\begin{equation}
\label{f3_37}
M_l M_r [ \mathbf{t}^{\mathbf{k}_j} ]_{\mathfrak{L}} = M_l [ \mathbf{t}^{\mathbf{k}_j + \vec e_r} ]_{\mathfrak{L}},\qquad
\forall j\in\Omega_0. 
\end{equation}
In general, it is not clear if the element $\mathbf{k}_j + \vec e_r =: \mathbf{k}_s$, with $s\in\{0,\ldots,\rho\}$,
has the property $s\in\Omega_l$. Thus, we can not apply
relation~(\ref{f3_7}) to get $[ \mathbf{t}^{\mathbf{k}_j + \vec e_r + \vec e_l} ]_{\mathfrak{L}}$. 
\textit{However, for an important type of admissible sets $K$, described in Example~\ref{e1_1}, part~2), this property
holds and we come to the commutativity of operators $M_k$.}
Then we can construct an atomic solution $\mu$ by relation~(\ref{f3_19}).

Assume that $s_{\mathbf{0}} \not= 0$. Then $H\not= \{ 0 \}$.
Apply the Gram-Schmidt orthogonalization process, removing linearly dependent elements, to the elements
$[ \mathbf{t}^{\mathbf{k}_j} ]_{\mathfrak{L}}$, $j\in\Omega_0$. Then we get an orthonormal basis
$\mathfrak{F} = \{ f_j \}_{j=0}^\tau$ in $H_0 = H$.

For an arbitrary admissible $K$, the commutativity of operators $M_j$ can be directly checked by using their matrices with respect to $\mathfrak{F}$.
We do not know, if the commutativity is true for all admissible sets $K$ (in the case of the dimensional stability).

Observe that conditions~(\ref{f2_25}),(\ref{f2_45}) can be verified by the standard tools. To verify the dimensional stability,
one can find projections of elements $[ \mathbf{t}^{\mathbf{k}_j} ]_{\mathfrak{L}}$, $j\in\{ 0,\ldots, r \}\backslash\Omega_0$,
on the subspace $H_0$, using an orthonormal basis in $H_0$. The norms of these projections should be zero.

\begin{example}
\label{e3_1}
Consider the truncated moment problem~(\ref{f1_1}) with $n=2$, 
$K = K_{2,2}$ (see~Example~\ref{e1_1}), $\mathcal{K} = K+K = K_{4,4}$, and the following moments:
$$ s_{(0,0)} = 3,\ s_{(0,1)} = s_{(0,2)} = s_{(0,3)} = s_{(0,4)} = 1, $$
$$ s_{(1,0)} = 4,\ s_{(1,1)} = s_{(1,2)} = s_{(1,3)} = s_{(1,4)} = 0, $$
$$ s_{(2,0)} = 4,\ s_{(2,1)} = s_{(2,2)} = s_{(2,3)} = s_{(2,4)} = 0, $$
$$ s_{(3,0)} = 16,\ s_{(3,1)} = s_{(3,2)} = s_{(3,3)} = s_{(3,4)} = 0, $$
$$ s_{(4,0)} = 32,\ s_{(4,1)} = s_{(4,2)} = s_{(4,3)} = s_{(4,4)} = 0. $$

Choose the following indexation in the set $K$:
$$ \mathbf{k}_0 = (0,0),\ \mathbf{k}_1 = (0,1),\ \mathbf{k}_2 = (0,2), $$   
$$ \mathbf{k}_3 = (1,0),\ \mathbf{k}_4 = (1,1),\ \mathbf{k}_5 = (1,2), $$   
$$ \mathbf{k}_6 = (2,0),\ \mathbf{k}_7 = (2,1),\ \mathbf{k}_8 = (2,2). $$   
Thus, we have $\rho = 8$. The matrix $\Gamma = (s_{\mathbf{k}_j + \mathbf{k}_m})_{m,j=0}^8$ has the following form:

\begin{equation}
\label{f3_45}
\Gamma =
\left(
\begin{array}{ccccccccc}
3 & 1 & 1 & 4 & 0 & 0 & 8 & 0 & 0 \\
1 & 1 & 1 & 0 & 0 & 0 & 0 & 0 & 0 \\
1 & 1 & 1 & 0 & 0 & 0 & 0 & 0 & 0 \\
4 & 0 & 0 & 8 & 0 & 0 & 16 & 0 & 0 \\
0 & 0 & 0 & 0 & 0 & 0 & 0 & 0 & 0 \\
0 & 0 & 0 & 0 & 0 & 0 & 0 & 0 & 0 \\
8 & 0 & 0 & 16 & 0 & 0 & 32 & 0 & 0 \\
0 & 0 & 0 & 0 & 0 & 0 & 0 & 0 & 0 \\
0 & 0 & 0 & 0 & 0 & 0 & 0 & 0 & 0 \end{array}
\right).
\end{equation}
The non-negativity of $\Gamma$ can be verified directly, by checking that the determinants of all submatrices, standing on the intersections
of rows and columns with the same indices, are non-negative.
The matrices $\Gamma_1, \Gamma_2, \widehat\Gamma_1, \widehat\Gamma_2$ have the following forms:
\begin{equation}
\label{f3_50}
\Gamma_1 =
\left(
\begin{array}{cccccc}
3 & 1 & 1 & 4 & 0 & 0  \\
1 & 1 & 1 & 0 & 0 & 0  \\
1 & 1 & 1 & 0 & 0 & 0  \\
4 & 0 & 0 & 8 & 0 & 0  \\
0 & 0 & 0 & 0 & 0 & 0  \\
0 & 0 & 0 & 0 & 0 & 0 \end{array}
\right),\quad
\Gamma_2 =
\left(
\begin{array}{cccccc}
3 & 1 & 4 & 0 & 8 & 0  \\
1 & 1 & 0 & 0 & 0 & 0  \\
4 & 0 & 8 & 0 & 16 & 0  \\
0 & 0 & 0 & 0 & 0 & 0  \\
8 & 0 & 16 & 0 & 32 & 0  \\
0 & 0 & 0 & 0 & 0 & 0\end{array}
\right),
\end{equation}

\begin{equation}
\label{f3_52}
\widehat\Gamma_1 =
\left(
\begin{array}{cccccc}
8 & 0 & 0 & 16 & 0 & 0 \\
0 & 0 & 0 & 0 & 0 & 0 \\
0 & 0 & 0 & 0 & 0 & 0 \\
16 & 0 & 0 & 32 & 0 & 0 \\
0 & 0 & 0 & 0 & 0 & 0 \\
0 & 0 & 0 & 0 & 0 & 0\end{array}
\right),\quad
\widehat\Gamma_2 =
\left(
\begin{array}{cccccc}
1 & 1 & 0 & 0 & 0 & 0 \\
1 & 1 & 0 & 0 & 0 & 0 \\
0 & 0 & 0 & 0 & 0 & 0 \\
0 & 0 & 0 & 0 & 0 & 0 \\
0 & 0 & 0 & 0 & 0 & 0 \\
0 & 0 & 0 & 0 & 0 & 0\end{array}
\right).
\end{equation}
The linear algebraic equation
$$ \Gamma_1 \vec x = 0,\quad \vec x = (x_1,\ldots,x_6)^T, $$
has the following solution:
$x_3,x_4,x_5,x_6$ are arbitrary complex numbers, $x_1 = -2 x_4$, $x_2 = - x_3 + 2 x_4$.
It is readily checked that for any solution it holds $\widehat \Gamma_1 \vec x = 0$.

\noindent
On the other hand, the linear algebraic equation
$$ \Gamma_2 \vec x = 0,\quad \vec x = (x_1,\ldots,x_6)^T, $$
has the following solution:
$x_2,x_4,x_5,x_6$ are arbitrary complex numbers, $x_1 = - x_2$, $x_3 = \frac{1}{2} x_2 - 2 x_5$.
It is readily checked that for any solution it holds $\widehat \Gamma_2 \vec x = 0$.

Thus, conditions~(\ref{f2_25}),(\ref{f2_45}) hold. Let us check the dimensional stability.
Consider the associated Hilbert space $H$.
For simplicity, we denote
$$ g_j = [ \mathbf{t}^{\mathbf{k}_j} ]_{\mathfrak{L}},\qquad j=0,\ldots, 8. $$
Observe that
$$ \Omega_0 = \{ 0, 1, 3, 4 \}. $$
Let us apply the Gram-Schmidt orthogonalization process, removing linearly dependent elements, to the sequence
$g_0, g_1, g_3, g_4$. Notice that all norms and scalar products are calculated by the moments:
\begin{equation}
\label{f3_54}
(g_j,g_r)_H = \left( 
[ \mathbf{t}^{\mathbf{k}_j} ]_{\mathfrak{L}},
[ \mathbf{t}^{\mathbf{k}_r} ]_{\mathfrak{L}}
\right)_H =
s_{\mathbf{k}_j + s_{\mathbf{k}_r}},\qquad j,r=0,1,\ldots, 8.
\end{equation}
We obtain an orthonormal basis 
$\mathfrak{F} = \{ f_0, f_1 \}$ in $H_0$, with
\begin{equation}
\label{f3_55}
f_0 = \frac{1}{\sqrt{3}} g_0,\quad
f_1 = \sqrt{\frac{3}{2}} \left( 
g_1 - \frac{1}{3} g_0
\right).
\end{equation}
Moreover, it turned out that 
\begin{equation}
\label{f3_57}
g_3 = 2 g_0 - 2 g_1,\quad g_4 = 0.
\end{equation}
It remains to verify that the projections of elements $g_2,g_5,g_6,g_7,g_8$ on $H_0$ coincide with the corresponding elements.
For example,
$$ g_2 - (g_2,f_0) f_0 - (g_2,f_1) f_1 = g_2 - g_1, $$
but
$$ \| g_2 - g_1 \|_H^2 = ( g_2 - g_1, g_2 - g_1 )_H = (g_2,g_2) - (g_2,g_1) - (g_1,g_2) + (g_1,g_1) = 0. $$ 
For other elements, we proceed in a similar way.
Consequently, the sequence $\mathcal{S} = (s_{\mathbf{k}})_{\mathbf{k}\in\mathcal{K}}$ is dimensionally stable.

Let us construct an atomic solution $\mu$ of the moment problem.
Observe that
$$ \Omega_1 = \{ 0, 1, 2, 3, 4, 5 \},\quad \Omega_2 = \{ 0, 1, 3, 4, 6, 7 \}. $$
The operators $M_1$ and $M_2$ act in the following way:
$$ M_1 g_0 = g_3 = 2g_0 - g_1,\quad M_1 g_1 = g_4 =0; $$
$$ M_2 g_0 = g_1,\quad M_2 g_1 = g_2 = g_1. $$
Therefore
$$ M_1 f_0 = \frac{4}{3} f_0 - \frac{2\sqrt{2}}{3} f_1,\quad M_1 f_1 = - \frac{2\sqrt{2}}{3} f_0 + \frac{2}{3} f_1; $$
$$ M_2 f_0 = \frac{1}{3} f_0 + \frac{\sqrt{2}}{3} f_1,\quad M_1 f_1 =  \frac{\sqrt{2}}{3} f_0 + \frac{2}{3} f_1. $$
The matrices $\mathcal{M}_1,\mathcal{M}_2$ of operators $M_1,M_2$, respectively, for the basis 
$\mathfrak{F}$ are:
$$ \mathcal{M}_1 = \left(
\begin{array}{cc} \frac{4}{3} & - \frac{2\sqrt{2}}{3}\\
- \frac{2\sqrt{2}}{3} & \frac{2}{3}\end{array}
\right),\quad
\mathcal{M}_2 = \left(
\begin{array}{cc} \frac{1}{3} &  \frac{\sqrt{2}}{3}\\
\frac{\sqrt{2}}{3} & \frac{2}{3}\end{array}
\right). $$
The matrix $\mathcal{M}_1$ has eigenvalues $\lambda_1 =  0, \lambda_2 = 2$, with eigenvectors respectively:
$$ \vec u_1 = \frac{1}{\sqrt{3}} (1,\sqrt{2})^T,\quad \vec u_2 = \frac{1}{\sqrt{3}} (-\sqrt{2}, 1)^T. $$
The matrix $\mathcal{M}_2$ has eigenvalues $\widetilde\lambda_1 =  0, \widetilde\lambda_2 = 1$, with eigenvectors respectively:
$$ \vec v_1 = \frac{1}{\sqrt{3}} (\sqrt{2},-1)^T,\quad \vec v_2 = \frac{1}{\sqrt{3}} (1,\sqrt{2})^T. $$
Denote
$$ \mathcal{H}_1 = \mathop{\rm Lin}\nolimits \left\{ \frac{1}{\sqrt{3}} (f_0 + \sqrt{2} f_1) \right\},\quad 
\mathcal{H}_2 = \mathop{\rm Lin}\nolimits \left\{ \frac{1}{\sqrt{3}} (-\sqrt{2} f_0 + f_1) \right\}, $$
$$ \widetilde{\mathcal{H}}_1 = \mathop{\rm Lin}\nolimits \left\{ \frac{1}{\sqrt{3}} (\sqrt{2} f_0 - f_1) \right\},\quad 
\widetilde{\mathcal{H}}_2 = \mathop{\rm Lin}\nolimits \left\{ \frac{1}{\sqrt{3}} (f_0 + \sqrt{2} f_1) \right\}. $$
Observe that the spectral measure $E(\delta)$ in relation~(\ref{f3_19}) can have jumps at points $(x,y)$ with $x\in \{ \lambda_1,\lambda_2 \}$,
$y\in\{ \widetilde{\lambda}_1, \widetilde{\lambda}_2 \}$.
The measure support is contained in this set of four points.
Thus, the measure $\mu$ has at most $4$ atoms. Notice that
$$ \mu(\{ (x,y) \}) = (E(\{ (x,y) \}) g_0, g_0)_H = (E_1(\{ x \}) E_2(\{ y \}) g_0, g_0)_H = $$
\begin{equation}
\label{f3_60}
= (E_2(\{ y \}) g_0, E_1(\{ x \}) g_0)_H,\qquad \forall (x,y)\in\mathbb{R}^2. 
\end{equation}
Observe that
$$ E_1(\{ 0 \}) g_0 = P_{\mathcal{H}_1} g_0 = \frac{1}{\sqrt{3}} (f_0 + \sqrt{2} f_1), $$
$$ E_1(\{ 2 \}) g_0 = P_{\mathcal{H}_2} g_0 = \frac{\sqrt{6}}{3} (\sqrt{2} f_0  - f_1), $$
$$ E_2(\{ 0 \}) g_0 = P_{\widetilde{\mathcal{H}}_1} g_0 = \frac{\sqrt{6}}{3} (\sqrt{2} f_0  - f_1), $$
$$ E_2(\{ 1 \}) g_0 = P_{\widetilde{\mathcal{H}}_2} g_0 = \frac{1}{\sqrt{3}} (f_0 + \sqrt{2} f_1). $$
By~(\ref{f3_60}) we conclude that the solution $\mu$ is $2$-atomic, having jumps $1$ and $2$  at points $(0,1)$
and $(2,0)$, respectively.

\end{example}

In the case, \textit{when the dimensional stability does not hold}, one can parametrize self-ajoint extensions $\widetilde M_j$ in a finite-dimensional
Hilbert space $\widetilde H\supseteq H$ for each $M_j$ separately ($j=1,\ldots,n$). Then one may study the commutativity 
of $M_j$s. For example, this can be done by the investigation of the commutativity of their (finite size) matrices 
with respect to an orthonormal basis in $\widetilde H$.

On the other hand, one can write conditions that ensure that all $M_j$, but one $M_{j_0}$, are selfadjoint. Then one can 
parametrize self-ajoint extensions $\widetilde M_{j_0}$ of $M_{j_0}$ inside $H$. Finally, it remains to check the commutativity of all operators,
using their matrices in an orthonormal basis. 

In the next section we shall rewrite the above ideas as a detailed algorithm for the case $n=2$.

\section{An algorithm for the truncated two-dimensional moment problem.}

In this section we shall describe an algorithm which, under certain conditions, allows to construct atomic solutions of
the moment problem~(\ref{f1_1}). For simplicity, we restrict ourselves to the case $n=2$.

\noindent
\textbf{Algorithm 1.}

\noindent
\textbf{The given data}: an admissible finite set $K\subset\mathbb{Z}^2_+$, 
$\mathcal{K} := K+K$ and a set of prescribed moments $\mathcal{S} = (s_{\mathbf{k}})_{\mathbf{k}\in\mathcal{K}}$.

\noindent
\textbf{Step 1.}
Choose and fix some indexation~(\ref{f2_9}) for the set $K$, with $\mathbf{k}_0 = \mathbf{0}$.

\noindent
\textbf{Step 2.}
Check conditions~(\ref{f2_25}) and (\ref{f2_45}). If they do not hold, then \textit{the moment problem
has no solutions} and we stop the algorithm.

\noindent
\textbf{Step 3.}
Consider the associated Hilbert space $H$, which is defined as in the paragraph following formula~(\ref{f2_25}).
Although this space consists of abstract elements (classes of the equivalence), all required numerical calculations will be performed 
by the following basic correlation property:
\begin{equation}
\label{f4_5}
\left( [ \mathbf{t}^{\mathbf{k}_j} ]_{\mathfrak{L}}, [ \mathbf{t}^{\mathbf{k}_m} ]_{\mathfrak{L}} \right)_H =
s_{\mathbf{k}_j + \mathbf{k}_m},\qquad j,m\in\{ 0,1,...,\rho \}.
\end{equation}

For $l=1,2$ we consider the multiplication operators $M_l$ as in~(\ref{f3_5}).
For convenience, we denote
\begin{equation}
\label{f4_6}
g_j = [ \mathbf{t}^{\mathbf{k}_j} ]_{\mathfrak{L}},\qquad j\in\{ 0,1,...,\rho \}. 
\end{equation}
Then
\begin{equation}
\label{f4_7}
\left( g_j, g_m \right)_H =
s_{\mathbf{k}_j + \mathbf{k}_m},\qquad j,m\in\{ 0,1,...,\rho \},
\end{equation}
and
\begin{equation}
\label{f4_9}
M_l g_j =
[ \mathbf{t}^{\mathbf{k}_j + \vec e_l} ]_{\mathfrak{L}} =: g_{\eta(l;j)},\qquad j\in\Omega_l;\quad l=1,2. 
\end{equation}

\noindent
\textbf{Step 4.} 
If $\| g_0 \|_H^2 = s_{\mathbf{k}_0 + \mathbf{k}_0} = s_\mathbf{0} = 0$, then the moment problem can not have any solution
different from $\mu=0$. In this case, if all the moments are zero then $\mu = 0$ is a solution, otherwise there are no solutions.
Thus, in the case $s_\mathbf{0} = 0$ we stop the algorithm.

\noindent
\textbf{Step 5.} (\textit{The construction of an orthonormal basis}).

Apply the Gram-Schmidt orthogonalization procedure to the sequence
$$ g_0, g_1, ..., g_\rho, $$
removing the linearly dependent elements, if they appear.
We get an orthonormal basis 
$$ \mathfrak{F} = \{ f_0, f_1, ..., f_{\rho'} \} $$
in the Hilbert space $H$ (where $0\leq \rho'\leq\rho$).
By the construction, an element $f_j$ is a linear combination of $g_k$s, with explicitly calculated coefficients.
Notice that $f_0 \not= 0$.

\noindent
\textbf{Step 6.} (\textit{The parametrization of extensions for each $M_l$ separately}).

Observe that $M_l$ is defined on elements $g_j$, $j\in \Omega_l$ ($l=1,2$). 
At first, define linear operators $\widetilde M_l$ on these elements in the same way.
Denote
\begin{equation}
\label{f4_15}
\Omega_l' := \{ 0,1,...,\rho \} \backslash \Omega_l,\qquad l=1,2. 
\end{equation}
For $l=1,2$ one should repeat the following procedure.

Choose an arbitrary element $g_k$, $k\in\Omega_l'$. Calculate the norm of its projection on $D(\widetilde M_l)$.
If $g_k\in D(\widetilde M_l)$ then we skip this element. Otherwise, we set
\begin{equation}
\label{f4_17}
\widetilde M_l g_k = \sum_{j=0}^{\rho'} (\alpha_{l;k,j} + \beta_{l;k,j} i) f_j,\qquad  \alpha_{l;k,j},\beta_{l;k,j}\in\mathbb{R}.
\end{equation}
Then we take another element $g_k$, $k\in\Omega_l'$, and proceed in a similar way. We continue this procedure to define
$\widetilde M_l$ for all $g_k$, $k\in\Omega_l'$.

Using the linearity, we construct some linear (but not necessary self-adjoint) extensions $\widetilde M_l$ of $M_l$ on the whole $H$ ($l=1,2$).
Notice that the case $D(M_l)=H$ was not excluded. The latter case  means that the corresponding parameters
$\alpha_{l;k,j},\beta_{l;k,j}$ are absent.

\noindent
\textbf{Step 7.} (\textit{The calculation of matrices of $\widetilde M_l$}).
Observe that each $f_j$ is a linear combination of $g_k$s (by the Gram-Schmidt orthogonalization):
\begin{equation}
\label{f4_19}
f_j = \sum_{k=0}^\rho c_{j;k} g_k,\qquad  c_{j;k}\in\mathbb{C};\ j\in\mathbb{Z}_{0,\rho'};
\end{equation}
and vice versa:
\begin{equation}
\label{f4_22}
g_j = \sum_{k=0}^{\rho'} d_{j;k} f_k,\qquad  d_{j;k}\in\mathbb{C};\ j\in\mathbb{Z}_{0,\rho}.
\end{equation}
Then
$$ \widetilde M_l f_j = \sum_{k=0}^\rho c_{j;k} \widetilde M_l g_k. $$ 
By~(\ref{f4_9}), (\ref{f4_17}) and~(\ref{f4_22}) we see that $\widetilde M_l f_j$
is a linear combination of $f_k$ with some coefficients, which may depend linearly on $\alpha_{l;k,j},\beta_{l;k,j}$.

\noindent
In the basis $\mathfrak{F}$, we calculate the matrices $\mathcal{M}_1$, $\mathcal{M}_2$ of $\widetilde M_1$ and $\widetilde M_2$, respectively. 
Thus, the coefficients of $\mathcal{M}_1, \mathcal{M}_2$ may depend 
\textit{linearly} on real parameters $\alpha_{l;k,j},\beta_{l;k,j}$.

\noindent
\textbf{Step 8.} (\textit{The check for the self-adjointness and the commutativity}).

The following conditions:
\begin{equation}
\label{f4_25}
\mathcal{M}_1 = \mathcal{M}_1^*,\quad \mathcal{M}_2 = \mathcal{M}_2^*,
\end{equation}
and
\begin{equation}
\label{f4_27}
\mathcal{M}_1 \mathcal{M}_2 = \mathcal{M}_2 \mathcal{M}_1, 
\end{equation}
ensure the self-adjointness and the commutativity of $\widetilde M_1$ and $\widetilde M_2$.

Conditions~(\ref{f4_25}) generate linear algebraic systems for unknown real parameters $\alpha_{l;k,j},\beta_{l;k,j}$. They
can be solved by elementary methods (e.g. by the Gauss elimination method).

\noindent
Substitute the general solutions of linear systems~(\ref{f4_25}) (which may depend on free real parameters) into relation~(\ref{f4_27}).
Observe that the coefficients of matrices $\mathcal{M}_1,\mathcal{M}_2$ may depend \textit{linearly} on these new free real parameters.
If parameters of $\mathcal{M}_1$ or parameters of $\mathcal{M}_2$ are absent, then we obtain a linear algebraic system of 
equations. In general, if we fix the parameters of $\mathcal{M}_2$, then we get a linear system with respect to parameters of $\mathcal{M}_2$.
Similarly, if we fix the parameters of $\mathcal{M}_1$, then we get a linear system with respect to parameters of $\mathcal{M}_2$.
If we can not get a solution on this way, we stop the algorithm.

Choose and fix arbitrary parameters $\alpha_{1;k,j},\beta_{1;k,j}$, $\alpha_{2;k,j},\beta_{2;k,j}$,
satisfying relations~(\ref{f4_25}),(\ref{f4_27}).
In what follows, we shall consider matrices $\mathcal{M}_1,\mathcal{M}_2$ corresponding to these parameters.

\noindent
\textbf{Step 9.}
Find all eigenvalues and eigenvectors of matrices $\mathcal{M}_1,\mathcal{M}_2$.

\noindent
\textbf{Step 10.}
Calculate the (atomic) solution of the moment problem by formula~(\ref{f3_19}).
Observe that the solution can have atoms at points $(x,y)$, where $x$ is an eigenvalue of $\mathcal{M}_1$,
$y$ is an eigenvalue of $\mathcal{M}_2$. Notice that
\begin{equation}
\label{f4_28}
\mu(\{ (x,y) \}) = (E(\{ (x,y) \}) g_0, g_0)_H = (E_2(\{ y \}) g_0, E_1(\{ x \}) g_0)_H,
\end{equation}
where $E_j(\delta)$ ($\delta\in\mathfrak{B}(\mathbb{R})$) is the spectral measure of the (bounded) self-adjoint operator $\widetilde M_j$,
$j=1,2$.

Let us illustrate this algorithm by the following examples.

\begin{example}
\label{e4_1}
Consider the truncated moment problem~(\ref{f1_1}) with $n=2$, 
$K = K_{1,1}$ (see~Example~\ref{e1_1}), $\mathcal{K} = K+K = K_{2,2}$, and the following moments:
$$ s_{(0,0)} = 4,\ s_{(0,1)} = 12,\ s_{(0,2)} = 48, $$
$$ s_{(1,0)} = 4,\ s_{(1,1)} = 12,\ s_{(1,2)} = 48, $$
$$ s_{(2,0)} = 4,\ s_{(2,1)} = 12,\ s_{(2,2)} = 48. $$

\noindent
\textbf{Step 1.}
Choose the following indexation in the set $K$:
$$ \mathbf{k}_0 = (0,0),\ \mathbf{k}_1 = (0,1),\ \mathbf{k}_2 = (1,0),\  \mathbf{k}_3 = (1,1). $$   
Thus, we have $\rho = 3$. 

\noindent
\textbf{Step 2.}
The matrix $\Gamma = (s_{\mathbf{k}_j + \mathbf{k}_m})_{m,j=0}^3$ has the following form:

\begin{equation}
\label{f4_30}
\Gamma =
\left(
\begin{array}{cccc}
4 & 12 & 4 & 12 \\
12 & 48 & 12 & 48 \\
4 & 12 & 4 & 12 \\
12 & 48 & 12 & 48 \end{array}
\right).
\end{equation}
The non-negativity of $\Gamma$ holds. It is verified by checking that the determinants of all submatrices, standing on the intersections
of rows and columns with the same indices, are non-negative.
Observe that
\begin{equation}
\label{f4_32}
\Omega_1 = \{ 0,1 \},\ \Omega_2 = \{ 0,2 \}. 
\end{equation}
The matrices $\Gamma_1, \Gamma_2, \widehat\Gamma_1, \widehat\Gamma_2$ have the following forms:
\begin{equation}
\label{f4_50}
\Gamma_1 =
\left(
\begin{array}{cc}
4 & 12 \\
12 & 48\end{array}
\right),\quad
\Gamma_2 =
\left(
\begin{array}{cc}
4 & 4 \\
4 & 4\end{array}
\right),
\end{equation}

\begin{equation}
\label{f4_52}
\widehat\Gamma_1 =
\left(
\begin{array}{cc}
4 & 12 \\
12 & 48 \end{array}
\right),\quad
\widehat\Gamma_2 =
\left(
\begin{array}{cc}
48 & 48 \\
48 & 48\end{array}
\right).
\end{equation}
Therefore conditions~(\ref{f2_45}) hold. 

\noindent
\textbf{Step 3.}
Consider the associated Hilbert space $H$, which is defined as in the paragraph following formula~(\ref{f2_25}).
Consider the multiplication operators $M_l$ as in~(\ref{f3_5}). Denote
$$ g_j = [ \mathbf{t}^{\mathbf{k}_j} ]_{\mathfrak{L}},\qquad j=0,1,2,3. $$ 
Notice that
\begin{equation}
\label{f4_55}
M_1 g_0 = g_2,\quad M_1 g_1 = g_3,
\end{equation}
\begin{equation}
\label{f4_57}
M_2 g_0 = g_1,\quad M_2 g_2 = g_3,
\end{equation}
and 
\begin{equation}
\label{f4_58}
D(M_1) = \mathop{\rm Lin}\nolimits \{ g_0, g_1 \},\quad
D(M_2) = \mathop{\rm Lin}\nolimits \{ g_0, g_2 \}.
\end{equation}

\noindent
\textbf{Step 4.}
In our case we have $\| g_0 \|_H^2 = s_{\mathbf{k}_0 + \mathbf{k}_0} = s_\mathbf{0} = 4\not= 0$.

\noindent
\textbf{Step 5.}
Let us apply the Gram-Schmidt orthogonalization process, removing linearly dependent elements, to the sequence
$g_0, g_1, g_2, g_3$. 
We shall use the property~(\ref{f4_7}).
We obtain that
\begin{equation}
\label{f4_61}
f_0 = \frac{1}{2} g_0,\quad f_1 = \frac{1}{2\sqrt{3}} \left(
g_1 - 3 g_0 \right),
\end{equation}
and
\begin{equation}
\label{f4_63}
g_2 = g_0,\quad g_3 = g_1.
\end{equation}
Therefore 
$\mathfrak{F} := \{ f_0, f_1 \}$ is an orthonormal basis in $H$, and $\rho' = 1$.

\noindent
\textbf{Step 6.}
Notice that
\begin{equation}
\label{f4_65}
\Omega_1' = \{ 2,3 \},\ \Omega_2' = \{ 1,3 \}. 
\end{equation}
By~(\ref{f4_58}) and (\ref{f4_63}) we see that $D(M_1) = H$ and $D(M_2) = \mathop{\rm Lin}\nolimits \{ g_0 \}$.
Therefore $\widetilde M_1 = M_1$.
Define $\widetilde M_2$ on $g_1$ in the following way:
\begin{equation}
\label{f4_67}
\widetilde M_2 g_1 = \sum_{j=0}^{1} (\alpha_{2;1,j} + \beta_{2;1,j} i) f_j,\qquad  \alpha_{2;1,j},\beta_{2;1,j}\in\mathbb{R}.
\end{equation}
Since $g_3 = g_1$, the procedure of the extension is finished. 

\noindent
\textbf{Step 7.}
By~(\ref{f4_55}),(\ref{f4_63}) we see that $M_1 = \widetilde M_1 = E_H$. Therefore, $\mathcal{M}_1$ is the identity matrix.
Let us calculate $\mathcal{M}_2$.
Observe that
\begin{equation}
\label{f4_71}
g_0 = 2 f_0,\quad g_1 = 6 f_0 + 2\sqrt{3} f_1. 
\end{equation}
By~(\ref{f4_61}),(\ref{f4_71}),(\ref{f4_57}),(\ref{f4_67}) we get
$$ \widetilde M_2 f_0 = 3 f_0 + \sqrt{3} f_1, $$
$$ \widetilde M_2 f_1 = \left(
\frac{1}{2\sqrt{3}} (\alpha_{2;1,0} + \beta_{2;1,0} i) - 18
\right) f_0
+  
\left(
\frac{1}{2\sqrt{3}} (\alpha_{2;1,1} + \beta_{2;1,1} i) - 6\sqrt{3}
\right)
f_1. $$
Then
\begin{equation}
\label{f4_73}
\mathcal{M}_2 =
\left(
\begin{array}{cc}
3 & \frac{1}{2\sqrt{3}} (\alpha_{2;1,0} + \beta_{2;1,0} i) - 18 \\
\sqrt{3} & \frac{1}{2\sqrt{3}} (\alpha_{2;1,1} + \beta_{2;1,1} i) - 6\sqrt{3} \end{array}
\right).
\end{equation}

\noindent
\textbf{Step 8.}
Conditions~(\ref{f4_25}) imply that
\begin{equation}
\label{f4_75}
\alpha_{2;1,0} = 36\sqrt{3} + 6,\quad \beta_{2;1,0} = \beta_{2;1,0} = 0. 
\end{equation}
Conditions~(\ref{f4_25}) are satisfied, since $\mathcal{M}_1$ is the identity matrix.
The real parameter $\alpha_{2;1,1}$ is free. We choose
$\alpha_{2;1,1} = 36 + 2\sqrt{3}$ to get
\begin{equation}
\label{f4_77}
\mathcal{M}_2 =
\left(
\begin{array}{cc}
3 & \sqrt{3} \\
\sqrt{3} & 1 \end{array}
\right).
\end{equation}

\noindent
\textbf{Step 9.}
The matrix $\mathcal{M}_1$ has an eigenvalue $\lambda_1 = 1$ and the eigensubspace
$\mathcal{H}_1 = H$.
The matrix $\mathcal{M}_2$ has eigenvalues $\widetilde\lambda_1 = 0$ and $\widetilde\lambda_2 = 4$, with eigensubspaces
$$ \widetilde{\mathcal{H}}_1 = \mathop{\rm Lin}\nolimits \left\{ -\frac{1}{2} f_0 + \frac{\sqrt{3}}{2} f_1 \right\},\quad
\widetilde{\mathcal{H}}_2 = \mathop{\rm Lin}\nolimits \left\{ \frac{\sqrt{3}}{2} f_0 + \frac{1}{2} f_1 \right\}, $$
respectively.

\noindent
\textbf{Step 10.}
Observe that
$$ E_2 (\{ 0 \}) g_0 = P_{\widetilde{\mathcal H}_1} g_0 = \frac{1}{2} f_0 - \frac{\sqrt{3}}{2} f_1, $$
$$ E_2 (\{ 4 \}) g_0 = P_{\widetilde{\mathcal H}_2} g_0 = \frac{3}{2} f_0 + \frac{\sqrt{3}}{2} f_1. $$
By formula~(\ref{f4_28}) we obtain that the solution $\mu$ is $2$-atomic with jumps $1$, $3$ at points
$(1,0)$ and $(1,4)$, respectively.

\end{example}

\begin{example}
\label{e4_2}
Consider the truncated moment problem~(\ref{f1_1}) with $n=2$, 
$K = K_{1,1}$ (see~Example~\ref{e1_1}), $\mathcal{K} = K+K = K_{2,2}$, and the following moments:
$$ s_{(0,0)} = 3,\ s_{(0,1)} = 2,\ s_{(0,2)} = 2, $$
$$ s_{(1,0)} = 3,\ s_{(1,1)} = 2,\ s_{(1,2)} = 2, $$
$$ s_{(2,0)} = 5,\ s_{(2,1)} = 4,\ s_{(2,2)} = 4. $$

\noindent
\textbf{Step 1.} The same as in the previous example.

\noindent
\textbf{Step 2.}
The matrix $\Gamma = (s_{\mathbf{k}_j + \mathbf{k}_m})_{m,j=0}^3$ has the following form:

\begin{equation}
\label{f4_80}
\Gamma =
\left(
\begin{array}{cccc}
3 & 2 & 3 & 2 \\
2 & 2 & 2 & 2 \\
3 & 2 & 5 & 4 \\
2 & 2 & 4 & 4 \end{array}
\right).
\end{equation}
The non-negativity of $\Gamma$ holds. 
Observe that
\begin{equation}
\label{f4_82}
\Omega_1 = \{ 0,1 \},\ \Omega_2 = \{ 0,2 \}. 
\end{equation}
The matrices $\Gamma_1, \Gamma_2, \widehat\Gamma_1, \widehat\Gamma_2$ have the following forms:
\begin{equation}
\label{f4_90}
\Gamma_1 =
\left(
\begin{array}{cc}
3 & 2 \\
2 & 2\end{array}
\right),\quad
\Gamma_2 =
\left(
\begin{array}{cc}
3 & 3 \\
3 & 5\end{array}
\right),
\end{equation}

\begin{equation}
\label{f4_92}
\widehat\Gamma_1 =
\left(
\begin{array}{cc}
5 & 4 \\
4 & 4 \end{array}
\right),\quad
\widehat\Gamma_2 =
\left(
\begin{array}{cc}
2 & 2 \\
2 & 4\end{array}
\right).
\end{equation}
Therefore conditions~(\ref{f2_45}) hold. 

\noindent
\textbf{Step 3.} The same as in the previous example.

\noindent
\textbf{Step 4.}
In our case we have $\| g_0 \|_H^2 = s_{\mathbf{k}_0 + \mathbf{k}_0} = s_\mathbf{0} = 3\not= 0$.

\noindent
\textbf{Step 5.}
Let us apply the Gram-Schmidt orthogonalization process, removing linearly dependent elements, to the sequence
$g_0, g_1, g_2, g_3$. 
We shall use the property~(\ref{f4_7}).
We obtain that
\begin{equation}
\label{f4_101}
f_0 = \frac{1}{\sqrt{3}} g_0,\quad f_1 = \sqrt{ \frac{3}{2} } \left(
g_1 - \frac{2}{3} g_0 \right),\quad f_2 = \frac{1}{\sqrt{2}} (g_2 - g_0),
\end{equation}
and
\begin{equation}
\label{f4_103}
g_3 = -g_0 + g_1 + g_2.
\end{equation}
Therefore 
$\mathfrak{F} := \{ f_0, f_1, f_2 \}$ is an orthonormal basis in $H$, and $\rho' = 2$.

\noindent
\textbf{Step 6.}
Notice that
\begin{equation}
\label{f4_105}
\Omega_1' = \{ 2,3 \},\ \Omega_2' = \{ 1,3 \}. 
\end{equation}

Define $\widetilde M_1$ on $g_2$ in the following way:
\begin{equation}
\label{f4_107}
\widetilde M_1 g_2 = \sum_{j=0}^{2} (\alpha_{1;2,j} + \beta_{1;2,j} i) f_j,\qquad  \alpha_{1;2,j},\beta_{1;2,j}\in\mathbb{R}.
\end{equation}

We define $\widetilde M_2$ on $g_1$ by the following formula:
\begin{equation}
\label{f4_109}
\widetilde M_2 g_1 = \sum_{j=0}^{2} (\alpha_{2;1,j} + \beta_{2;1,j} i) f_j,\qquad  \alpha_{2;1,j},\beta_{2;1,j}\in\mathbb{R}.
\end{equation}
Since $g_3 = -g_0 + g_1 + g_2$, the procedure of the extension is finished. 
Notice that we have $12$ free real parameters at this moment.

\noindent
\textbf{Step 7.}
Observe that
$$ g_0 = \sqrt{3} f_0,\quad g_1 = \frac{2}{\sqrt{3}} f_0 + \sqrt{ \frac{2}{3} } f_1,\quad g_2 = \sqrt{3} f_0 + \sqrt{2} f_2,  $$
\begin{equation}
\label{f4_111}
g_3 = \frac{2}{\sqrt{3}} f_0 + \sqrt{ \frac{2}{3} } f_1 + \sqrt{2} f_2. 
\end{equation}
By~(\ref{f4_101}),(\ref{f4_111}),(\ref{f4_107}),(\ref{f4_109}),(\ref{f4_55}) and (\ref{f4_57}) we get
$$ \widetilde M_1 f_0 = f_0 + \sqrt{ \frac{2}{3} } f_2,\quad \widetilde M_1 f_1 = f_1 + \frac{1}{ \sqrt{3} } f_2, $$
$$ \widetilde M_1 f_2 = 
\frac{1}{ \sqrt{2} } 
\left(
\alpha_{1;2,0} + \beta_{1;2,0} i - \sqrt{3}
\right) f_0
+  
\frac{1}{ \sqrt{2} } 
\left(
\alpha_{1;2,1} + \beta_{1;2,1} i
\right)
f_1
+ $$
$$ + 
\frac{1}{ \sqrt{2} } 
\left(
\alpha_{1;2,2} + \beta_{1;2,2} i - \sqrt{2}
\right) f_2; $$

$$ \widetilde M_2 f_0 = \frac{2}{3} f_0 + \frac{ \sqrt{2} }{3} f_1, $$
$$ \widetilde M_2 f_1 = 
\left(
\sqrt{ \frac{3}{2} } (\alpha_{2;1,0} + \beta_{2;1,0} i) - \frac{ 2\sqrt{2} }{3}
\right) f_0
+  $$
$$ + \left(
\sqrt{ \frac{3}{2} } (\alpha_{2;1,1} + \beta_{2;1,1} i) - \frac{ 2 }{3}
\right) f_1
+
\sqrt{ \frac{3}{2} }
\left(
\alpha_{2;1,2} + \beta_{2;1,2} i
\right)
f_2, $$
$$ \widetilde M_2 f_2 = f_2. $$ 
Then
\begin{equation}
\label{f4_114}
\mathcal{M}_1 =
\left(
\begin{array}{ccc}
1 & 0 & \frac{1}{ \sqrt{2} } 
\left(
\alpha_{1;2,0} + \beta_{1;2,0} i - \sqrt{3}
\right) \\
0 & 1 & \frac{1}{ \sqrt{2} } 
\left(
\alpha_{1;2,1} + \beta_{1;2,1} i
\right) \\
\sqrt{\frac{2}{3}} & \frac{1}{\sqrt{3}} & \frac{1}{ \sqrt{2} } 
\left(
\alpha_{1;2,2} + \beta_{1;2,2} i - \sqrt{2}
\right) \end{array}
\right),
\end{equation}
\begin{equation}
\label{f4_116}
\mathcal{M}_2 =
\left(
\begin{array}{ccc}
\frac{2}{3} & 
\sqrt{ \frac{3}{2} } (\alpha_{2;1,0} + \beta_{2;1,0} i) - \frac{ 2\sqrt{2} }{3} & 0 \\
\frac{\sqrt{2}}{3} & 
\sqrt{ \frac{3}{2} } (\alpha_{2;1,1} + \beta_{2;1,1} i) - \frac{ 2 }{3} & 0\\
0 & \sqrt{ \frac{3}{2} }
\left(
\alpha_{2;1,2} + \beta_{2;1,2} i
\right) & 1\end{array}\right).
\end{equation}

\noindent
\textbf{Step 8.}
Conditions~(\ref{f4_25}) imply that
$$ \beta_{1;2,j} = \beta_{2;1,j} = 0,\qquad j=0,1,2; $$
\begin{equation}
\label{f4_118}
\alpha_{1;2,0} = \frac{5}{ \sqrt{3} },\ \alpha_{1;2,1} = \sqrt{ \frac{2}{3} },\ 
\alpha_{2;1,0} = \frac{2}{ \sqrt{3} },\ \alpha_{2;1,2} = 0. 
\end{equation}
It remains two free real parameters: $\alpha_{1;2,2}$ and $\alpha_{2;1,1}$.
Matrices $\mathcal{M}_1$, $\mathcal{M}_2$ take the following form:
\begin{equation}
\label{f4_122}
\mathcal{M}_1 =
\left(
\begin{array}{ccc}
1 & 0 & \sqrt{\frac{2}{3}} \\
0 & 1 & \frac{1}{ \sqrt{3} } \\
\sqrt{\frac{2}{3}} & \frac{1}{\sqrt{3}} & 
\frac{1}{ \sqrt{2} } 
\alpha_{1;2,2} - 1 \end{array}
\right),
\end{equation}
\begin{equation}
\label{f4_124}
\mathcal{M}_2 =
\left(
\begin{array}{ccc}
\frac{2}{3} & 
\frac{ \sqrt{2} }{3} & 0 \\
\frac{\sqrt{2}}{3} & 
\sqrt{ \frac{3}{2} } \alpha_{2;1,1} - \frac{ 2 }{3} & 0\\
0 & 0 & 1\end{array}\right).
\end{equation}
Condition~(\ref{f4_27}) will be satisfied if 
$\alpha_{2;1,1} = \sqrt{ \frac{2}{3} }$. 
It remains one free real parameter $\alpha_{1;2,2}$. We set
$\alpha_{1;2,2} = 2\sqrt{2}$.
Therefore
\begin{equation}
\label{f4_126}
\mathcal{M}_1 =
\left(
\begin{array}{ccc}
1 & 0 & \sqrt{\frac{2}{3}} \\
0 & 1 & \frac{1}{ \sqrt{3} } \\
\sqrt{\frac{2}{3}} & \frac{1}{\sqrt{3}} & 1 \end{array}
\right),\quad
\mathcal{M}_2 =
\left(
\begin{array}{ccc}
\frac{2}{3} & 
\frac{ \sqrt{2} }{3} & 0 \\
\frac{\sqrt{2}}{3} & 
\frac{ 1 }{3} & 0\\
0 & 0 & 1\end{array}\right).
\end{equation}

\noindent
\textbf{Step 9.}
The matrix $\mathcal{M}_1$ has eigenvalues $\lambda_0 = 0$, $\lambda_1 = 1$ and $\lambda_2 = 2$, with eigensubspaces
$$ \mathcal{H}_0 = \mathop{\rm Lin}\nolimits \left\{ -\frac{1}{\sqrt{3}} f_0 - \frac{1}{ \sqrt{6} } f_1 
+ \frac{1}{ \sqrt{2} } f_2 \right\},\quad
\mathcal{H}_1 = \mathop{\rm Lin}\nolimits \left\{ \frac{1}{\sqrt{3}} f_0 - \sqrt{ \frac{2}{3} } f_1 \right\}, $$
$$ \mathcal{H}_2 = \mathop{\rm Lin}\nolimits \left\{ \frac{1}{\sqrt{3}} f_0 + \frac{1}{ \sqrt{6} } f_1 
+ \frac{1}{ \sqrt{2} } f_2 \right\}, $$
respectively.
The matrix $\mathcal{M}_2$ has eigenvalues $\widetilde\lambda_0 = 0$ and $\widetilde\lambda_1 = 1$, with eigensubspaces
$$ \widetilde{\mathcal{H}}_0 = \mathop{\rm Lin}\nolimits \left\{ -\frac{1}{ \sqrt{3} } f_0 + \sqrt{ \frac{2}{3} } f_1 \right\}, $$
$$ \widetilde{\mathcal{H}}_1 = 
\mathop{\rm Lin}\nolimits \left\{ \frac{1}{\sqrt{2}} f_0 + \frac{1}{2} f_1 + \frac{1}{2} f_2;\
-\frac{1}{ \sqrt{6} } f_0 - \frac{1}{ 2\sqrt{3} } f_1 
+ \frac{ \sqrt{3} }{2} f_2 \right\}, $$
respectively.

\noindent
\textbf{Step 10.}
Observe that
$$ E_1 (\{ 0 \}) g_0 = P_{\mathcal H_0} g_0 = \frac{1}{ \sqrt{3} } f_0 - \frac{1}{ \sqrt{6} } f_1 - \frac{1}{ \sqrt{2} } f_2, $$
$$ E_1 (\{ 1 \}) g_0 = P_{\mathcal H_1} g_0 = \frac{1}{ \sqrt{3} } f_0 - \sqrt{ \frac{2}{3} } f_1, $$
$$ E_1 (\{ 2 \}) g_0 = P_{\mathcal H_2} g_0 = \frac{1}{ \sqrt{3} } f_0 + \frac{1}{ \sqrt{6} } f_1 + \frac{1}{ \sqrt{2} } f_2; $$
$$ E_2 (\{ 0 \}) g_0 = P_{\widetilde{\mathcal H}_0} g_0 = \frac{1}{ \sqrt{3} } f_0 - \sqrt{ \frac{2}{3} } f_1, $$
$$ E_2 (\{ 1 \}) g_0 = P_{\widetilde{\mathcal H}_1} g_0 = \frac{ 2\sqrt{3} }{3} f_0 + \frac{ \sqrt{6} }{3} f_1. $$
By formula~(\ref{f4_28}) we obtain that the solution $\mu$ is $3$-atomic with unit jumps at points
$(0,1)$, $(1,0)$ and $(2,1)$.

\end{example}

\begin{remark}
\label{r4_1}
Observe that in Algorithm~1 we restricted ourselves by considering possible extensions $\widetilde M_j$ of $M_j$ inside the original Hilbert
space $H$. Instead of $H$ one can consider any finite-dimensional Hilbert space $\widetilde H\supseteq H$, and construct  
possible extensions $\widetilde M_j$ of $M_j$ in $\widetilde H$.

\end{remark}

\noindent
\textbf{Acknowledgements.} The author is grateful to Prof. Vasilescu for a useful discussion on the moment problems.


\begin{center}
{\large\bf 
The operator approach to the truncated multidimensional moment problem.}
\end{center}
\begin{center}
{\bf S.M. Zagorodnyuk}
\end{center}

We study the truncated multidimensional moment problem with a general type of truncations. The operator approach to
the moment problem is presented. A way to construct atomic solutions of the moment problem is indicated.

}

\noindent
Address:

V. N. Karazin Kharkiv National University \newline\indent
School of Mathematics and Computer Sciences \newline\indent
Department of Higher Mathematics and Informatics \newline\indent
Svobody Square 4, 61022, Kharkiv, Ukraine

Sergey.M.Zagorodnyuk@gmail.com; Sergey.M.Zagorodnyuk@univer.kharkov.ua

\end{document}